*Research Article*

# Startpoints and $(\alpha, \gamma)$-Contractions in Quasi-Pseudometric Spaces

### Yaé Ulrich Gaba

*Department of Mathematics and Applied Mathematics, University of Cape Town, Rondebosch 7701, South Africa*

Correspondence should be addressed to Yaé Ulrich Gaba; gabayae2@gmail.com

We introduce the concept of *startpoint* and *endpoint* for multivalued maps defined on a quasi-pseudometric space. We investigate the relation between these new concepts and the existence of fixed points for these set valued maps.

*Dedicated to my beloved Clémence on the occasion of her 25th birthday*

## 1. Introduction

In the last few years there has been a growing interest in the theory of quasi-metric spaces and other related structures such as quasi-normed cones and asymmetric normed linear spaces (see, e.g., [1]), because such a theory provides an important tool and a convenient framework in the study of several problems in theoretical computer science, applied physics, approximation theory, and convex analysis. Many works on general topology have been done in order to extend the well-known results of the classical theory. In particular, various types of completeness are studied in [2], showing, for instance, that the classical concept of *Cauchy sequences* can be accordingly modified. In the same reference, which uses an approach based on uniformities, the *bicompletion* of a $T_0$-quasi-pseudometric has been explored. It is worth mentioning that, in the fixed point theory, *completeness* is a key element, since most of the constructed sequences will be assumed to have a *Cauchy type* property.

It is the aim of this paper to continue the study of quasi-pseudometric spaces by proving some fixed point results and investigating a bit more the behaviour of set-valued mappings. Thus, in Section 3 a suitable notion of $(\alpha, \gamma)$-contractive mapping is given for self-mappings defined on quasi-pseudometric spaces and some fixed point results are discussed. In Sections 4 and 5, the notions of *startpoint* and

*endpoint* for set-valued mappings are introduced and different variants of such concepts, as well as their connections with the fixed point of a multivalued map, are characterized.

For recent results in the theory of asymmetric spaces, the reader is referred to [3–8].

## 2. Preliminaries

*Definition 1.* Let $X$ be a nonempty set. A function $d : X \times X \to [0, \infty)$ is called a *quasi-pseudometric* on $X$ if

   (i) $d(x, x) = 0$ for all $x \in X$,

   (ii) $d(x, z) \leq d(x, y) + d(y, z)$ for all $x, y, z \in X$.

Moreover, if $d(x, y) = 0 = d(y, x) \Rightarrow x = y$, then $d$ is said to be a $T_0$-*quasi-pseudometric*. The latter condition is referred to as the $T_0$-condition.

*Remark 2.* (i) Let $d$ be a quasi-pseudometric on $X$; then the map $d^{-1}$ defined by $d^{-1}(x, y) = d(y, x)$ whenever $x, y \in X$ is also a quasi-pseudometric on $X$, called the *conjugate* of $d$. In the literature, $d^{-1}$ is also denoted by $d^t$ or $\bar{d}$.

(ii) It is easy to verify that the function $d^s$ defined by $d^s := d \vee d^{-1}$, that is, $d^s(x, y) = \max\{d(x, y), d(y, x)\}$, defines a metric on $X$ whenever $d$ is a $T_0$-quasi-pseudometric on $X$.

Let $(X, d)$ be a quasi-pseudometric space. For $x \in X$ and $\varepsilon > 0$,

$$B_d (x, \varepsilon) = \{ y \in X : d(x, y) < \varepsilon \} \tag{1}$$

denotes the open $\varepsilon$-ball at $x$. The collection of all such balls yields a base for the topology $\tau(d)$ induced by $d$ on $X$. Hence, for any $A \in X$, we will, respectively, denote by $\text{int}_{\tau(d)} A$ and $\text{cl}_{\tau(d)} A$ the interior and the closure of the set $A$ with respect to the topology $\tau(d)$.

Similarly, for $x \in X$ and $\varepsilon \geq 0$,

$$C_d (x, \varepsilon) = \{ y \in X : d(x, y) \leq \varepsilon \} \tag{2}$$

denotes the closed $\varepsilon$-ball at $x$. We will say that a subset $E \subset X$ is *join-closed* if it is $\tau(d^s)$-closed, that is, closed with respect to the topology generated by $d^s$. The topology $\tau(d^s)$ is finer than the topologies $\tau(d)$ and $\tau(d^{-1})$.

*Definition 3.* Let $(X, d)$ be a quasi-pseudometric space. The convergence of a sequence $(x_n)$ to $x$ with respect to $\tau(d)$, called *$d$-convergence* or *left-convergence* and denoted by $x_n \xrightarrow{d} x$, is defined in the following way:

$$x_n \xrightarrow{d} x \iff d(x, x_n) \longrightarrow 0. \tag{3}$$

Similarly, the convergence of a sequence $(x_n)$ to $x$ with respect to $\tau(d^{-1})$, called *$d^{-1}$-convergence* or *right-convergence* and denoted by $x_n \xrightarrow{d^{-1}} x$, is defined in the following way:

$$x_n \xrightarrow{d^{-1}} x \iff d(x_n, x) \longrightarrow 0. \tag{4}$$

Finally, in a quasi-pseudometric space $(X, d)$, we will say that a sequence $(x_n)$ *$d^s$-converges* to $x$ if it is both left and right convergent to $x$, and we denote it by $x_n \xrightarrow{d^s} x$ or $x_n \to x$ when there is no confusion. Hence,

$$x_n \xrightarrow{d^s} x \iff x_n \xrightarrow{d} x, \qquad x_n \xrightarrow{d^{-1}} x. \tag{5}$$

*Definition 4.* A sequence $(x_n)$ in a quasi-pseudometric $(X, d)$ is called

(a) *left $d$-Cauchy* if for every $\epsilon > 0$, there exist $x \in X$ and $n_0 \in \mathbb{N}$ such that

$$\forall n \geq n_0 \quad d(x, x_n) < \epsilon; \tag{6}$$

(b) *left $K$-Cauchy* if for every $\epsilon > 0$, there exists $n_0 \in \mathbb{N}$ such that

$$\forall n, k : n_0 \leq k \leq n \quad d(x_k, x_n) < \epsilon; \tag{7}$$

(c) *$d^s$-Cauchy* if for every $\epsilon > 0$, there exists $n_0 \in \mathbb{N}$ such that

$$\forall n, k \geq n_0 \quad d(x_n, x_k) < \epsilon. \tag{8}$$

Dually, we define, in the same way, *right $d$-Cauchy* and *right $K$-Cauchy* sequences.

*Remark 5.* Consider the following:

(i) $d^s$-Cauchy $\Rightarrow$ left $K$-Cauchy $\Rightarrow$ left $d$-Cauchy. The same implications hold for the corresponding right notions. None of the above implications is reversible.

(ii) A sequence is left $d$-Cauchy with respect to $d$ if and only if it is right $K$-Cauchy with respect to $d^{-1}$.

(iii) A sequence is left $K$-Cauchy with respect to $d$ if and only if it is right $K$-Cauchy with respect to $d^{-1}$.

(iv) A sequence is $d^s$-Cauchy if and only if it is both left and right $K$-Cauchy.

*Definition 6.* A quasi-pseudometric space $(X, d)$ is called

(i) *left-$K$-complete* provided that any left $K$-Cauchy sequence is $d$-convergent,

(ii) *left Smyth sequentially complete* if any left $K$-Cauchy sequence is $d^s$-convergent.

The dual notions of *right-completeness* are easily derived from the above definition.

*Definition 7.* A $T_0$-quasi-pseudometric space $(X, d)$ is called *bicomplete* provided that the metric $d^s$ on $X$ is complete.

As usual, a subset $A$ of a quasi-pseudometric space $(X, d)$ will be called *bounded* provided that there exists a positive real constant $M$ such that $d(x, y) < M$ whenever $x, y \in A$. This is equivalent to saying that there exist $x_0 \in X$, $r, s \geq 0$ such that $A \subseteq C_d(x_0, r) \cap C_{d^{-1}}(x_0, s)$.

We also define the *diameter* $\delta(A)$ of $A$ by $\delta(A) := \sup\{d(x, y) : x, y \in A\}$. Hence, $A$ is bounded if and only if $\delta(A) < \infty$. It is not difficult to see that this definition coincides with that of a bounded set in a metric space.

Let $(X, d)$ be a quasi-pseudometric space. We set $\mathscr{P}_0(X) := 2^X \setminus \{0\}$ where $2^X$ denotes the power set of $X$. For $x \in X$ and $A, B \in \mathscr{P}_0(X)$, we define

$$\begin{aligned} d(x, A) &= \inf \{ d(x, a), \ a \in A \}, \\ d(A, x) &= \inf \{ d(a, x), \ a \in A \}, \end{aligned} \tag{9}$$

and we define $H(A, B)$ by

$$H(A, B) = \max \left\{ \sup_{a \in A} d(a, B), \sup_{b \in B} d(A, b) \right\}. \tag{10}$$

Then $H$ is an extended quasi-pseudometric on $\mathscr{P}_0(X)$. Moreover, we know from [9] that the restriction of $H$ to $S_{\text{cl}}(X) = \{ A \subseteq X : A = (\text{cl}_{\tau(d)} A) \cap (\text{cl}_{\tau(d^{-1})} A) \}$ is an extended $T_0$-quasi-pseudometric. We will denote by $CB(X)$ the collection of all nonempty bounded and $\tau(d^s)$-closed subsets of $X$.

We complete this section by the following lemma.

**Lemma 8.** Let $(X, d)$ be a quasi-pseudometric space. For every fixed $x \in X$, the mapping $y \mapsto d(x, y)$ is $\tau(d)$-upper semicontinuous ($\tau(d)$-usc in short) and $\tau(d^{-1})$-lower semicontinuous ($\tau(d^{-1})$-lsc in short). For every fixed $y \in X$, the mapping $x \mapsto d(x, y)$ is $\tau(d)$-lsc and $\tau(d^{-1})$-usc.

*Proof.* To prove that $d(x, \cdot)$ is $\tau(d)$-usc and $\tau(d^{-1})$-lsc, we have to show that the set $\{y \in X : d(x, y) < \alpha\}$ is $\tau(d)$-open and $\{y \in X : d(x, y) > \alpha\}$ is $\tau(d^{-1})$-open, for every $\alpha \in \mathbb{R}$, properties that are easy to check.

Indeed, for $y \in X$ such that $d(x, y) < \alpha$, let $r := \alpha - d(x, y) > 0$. If $z \in X$ is such that $d(y, z) < r$, then

$$d(x, z) \leq d(x, y) + (y, z) < d(x, y) + r = \alpha, \qquad (11)$$

showing that $B_d(y, r) \subset \{y \in X : d(x, y) < \alpha\}$.

Similarly, for $y \in X$ with $d(x, y) > \alpha$ take $r := d(x, y) - \alpha > 0$. If $z \in X$ satisfies $d(z, y) = d^{-1}(y, z) < r$, then

$$d(x, y) \leq d(x, z) + (z, y) < d(x, z) + r, \qquad (12)$$

so that $d(x, z) > d(x, y) - r = \alpha$. Consequently, $B_{d^{-1}}(y, r) \subset \{y \in X : d(x, y) > \alpha\}$.                                  $\square$

## 3. Some First Results

We begin by recalling the following.

*Definition 9.* A function $\gamma : [0, \infty) \to [0, \infty)$ is called a $(c)$-comparison function if

($\gamma_1$) $\gamma$ is nondecreasing;

($\gamma_2$) $\sum_{n=1}^{\infty} \gamma^n(t) < \infty$ for all $t > 0$, where $\gamma^n$ is the $n$th iterate of $\gamma$.

We will denote by $\Gamma$ the set of such functions. Note that for any $\gamma \in \Gamma$, $\gamma(t) < t$ for any $t > 0$.

We then introduce the following definitions.

*Definition 10.* Let $(X, d)$ be a quasi-pseudometric type space. A function $T : X \to X$ is called *$d$-sequentially continuous* or *left-sequentially continuous* if for any $d$-convergent sequence $(x_n)$ with $x_n \xrightarrow{d} x$, the sequence $(Tx_n)$ $d$-converges to $Tx$; that is, $Tx_n \xrightarrow{d} Tx$.

Similarly, we define a *$d^{-1}$-sequentially continuous* or *right-sequentially continuous* function.

*Definition 11.* Let $(X, d)$ be a quasi-pseudometric space, and let $T : X \to X$ and $\alpha : X \times X \to [0, \infty)$ be mappings. We say that $T$ is $\alpha$-admissible if

$$\alpha(x, y) \geq 1 \Longrightarrow \alpha(Tx, Ty) \geq 1, \qquad (13)$$

whenever $x, y \in X$.

*Definition 12.* Let $(X, d)$ be a quasi-pseudometric space and let $T : X \to X$ be a mapping. We say that $T$ is an $(\alpha, \gamma)$-contractive mapping if there exist two functions $\alpha : X \times X \to [0, \infty)$ and $\gamma \in \Gamma$ such that

$$\alpha(x, y) d(Tx, Ty) \leq \gamma(d(x, y)), \qquad (14)$$

whenever $x, y \in X$.

We now state the first fixed point theorem.

**Theorem 13.** *Let $(X, d)$ be a Hausdorff left $K$-complete $T_0$-quasi-pseudometric space. Suppose that $T : X \to X$ is an $(\alpha, \gamma)$-contractive mapping which satisfies the following conditions:*

(i) *$T$ is $\alpha$-admissible;*

(ii) *there exists $x_0 \in X$ such that $\alpha(x_0, Tx_0) \geq 1$;*

(iii) *$T$ is $d$-sequentially continuous.*

*Then $T$ has a fixed point.*

*Proof.* By (ii), there exists $x_0 \in X$ such that $\alpha(x_0, Tx_0) \geq 1$. Let us define the sequence $(x_n)$ in $X$ by $x_{n+1} = Tx_n$ for all $n = 0, 1, 2, \ldots$. Without loss of generality, we can always assume that $x_n \neq x_{n+1}$ for all $n \in \mathbb{N}$, since if $x_{n_0} = x_{n_0+1}$ for some $n_0 \in \mathbb{N}$, the proof is complete.                      $\square$

From assumption (i), we derive

$$\alpha(x_0, Tx_0) = \alpha(x_0, x_1) \geq 1$$
$$\Longrightarrow \alpha(Tx_0, Tx_1) = \alpha(x_1, x_2) \geq 1. \qquad (15)$$

Recursively, we get

$$\alpha(x_n, x_{n+1}) \geq 1 \quad \forall n = 0, 1, 2, \ldots. \qquad (16)$$

Since $T$ is $(\alpha, \gamma)$-contractive, we can write

$$d(x_n, x_{n+1}) = d(Tx_{n-1}, Tx_n)$$
$$\leq \alpha(x_{n-1}, x_n) d(Tx_{n-1}, Tx_n) \qquad (17)$$
$$\leq \gamma(d(x_{n-1}, x_n)),$$

for all $n = 0, 1, 2, \ldots$. Inductively, we obtain

$$d(x_n, x_{n+1}) \leq \gamma^n(d(x_0, x_1)), \quad n = 1, 2, \ldots. \qquad (18)$$

Therefore, for any $k \geq 1$, using the triangle inequality, we get

$$d(x_n, x_{n+k}) \leq d(x_n, x_{n+1}) + d(x_{n+1}, x_{n+2})$$
$$+ \cdots + d(x_{n+k-1}, x_{n+k})$$
$$\leq \sum_{i=n}^{n+k-1} \gamma^i(d(x_0, x_1)) \qquad (19)$$
$$\leq \sum_{i=n}^{\infty} \gamma^i(d(x_0, x_1)).$$

Letting $n \to \infty$, we derive that $d(x_n, x_{n+k}) \to 0$. Hence, $(x_n)$ is a left $K$-Cauchy sequence. Since $(X, d)$ is left $K$-complete and $T$ $d$-sequentially continuous, there exists $x^*$ such that $x_n \xrightarrow{d} x^*$ and $x_{n+1} \xrightarrow{d} Tx^*$. Since $X$ is Hausdorff, we have that $x^* = Tx^*$.

**Corollary 14.** *Let $(X, d)$ be a Hausdorff right $K$-complete $T_0$-quasi-pseudometric space. Suppose that $T : X \rightarrow X$ is an $(\alpha, \gamma)$-contractive mapping which satisfies the following conditions:*

(i) *$T$ is $\alpha$-admissible;*

(ii) *there exists $x_0 \in X$ such that $\alpha(Tx_0, x_0) \geq 1$;*

(iii) *$T$ is $d^{-1}$-sequentially continuous.*

*Then $T$ has a fixed point.*

**Corollary 15.** *Let $(X, d)$ be a bicomplete quasi-pseudometric space. Suppose that $T : X \rightarrow X$ is an $(\alpha, \gamma)$-contractive mapping which satisfies the following conditions:*

(i) *$T$ is $\alpha$-admissible and the function $\alpha$ is symmetric; that is, $\alpha(a, b) = \alpha(b, a)$ for any $a, b \in X$;*

(ii) *there exists $x_0 \in X$ such that $\alpha(Tx_0, x_0) \geq 1$;*

(iii) *$T$ is $d^s$-sequentially continuous.*

*Then $T$ has a fixed point.*

*Proof.* Following the proof of Theorem 13, it is clear that the sequence $(x_n)$ in $X$ defined by $x_{n+1} = Tx_n$ for all $n = 0, 1, 2, \ldots$ is $d^s$-Cauchy. Since $(X, d)$ is complete and $T$ sequentially continuous, there exists $x^*$ such that $x_n \xrightarrow{d^s} x^*$ and $x_{n+1} \xrightarrow{d^s} Tx^*$. Since $(X, d^s)$ is Hausdorff, we have that $x^* = Tx^*$. □

*Remark 16.* In fact, we do not need $\alpha$ to be symmetric. It is enough, for the result to be true, to have a point $x_0 \in X$ for which $\alpha(x_0, Tx_0) \geq 1$ and $\alpha(Tx_0, x_0) \geq 1$.

We conclude this section by the following results which are in fact consequences of Theorem 13.

**Theorem 17.** *Let $(X, d)$ be a Hausdorff left $K$-complete $T_0$-quasi-pseudometric space. Suppose that $T : X \rightarrow X$ is an $(\alpha, \gamma)$-contractive mapping which satisfies the following conditions:*

(i) *$T$ is $\alpha$-admissible;*

(ii) *there exists $x_0 \in X$ such that $\alpha(x_0, Tx_0) \geq 1$;*

(iii) *if $(x_n)$ is a sequence in $X$ such that $\alpha(x_n, x_{n+1}) \geq 1$ for all $n = 1, 2, \ldots$ and $x_n \xrightarrow{d} x$, then there exists a subsequence $(x_{n(k)})$ of $(x_n)$ such that $\alpha(x_{n(k)}, x) \geq 1$ for all $k$.*

*Then $T$ has a fixed point.*

*Proof.* Following the proof of Theorem 13, we know that the sequence $(x_n)$ defined by $x_{n+1} = Tx_n$ for all $n = 0, 1, 2, \ldots$ $d$-converges to some $x^*$ and satisfies $\alpha(x_n, x_{n+1}) \geq 1$ for $n \geq 1$. From condition (iii), we know that there exists a subsequence $(x_{n(k)})$ of $(x_n)$ such that $\alpha(x_{n(k)}, x^*) \geq 1$ for all $k$. Since $T$ is an $(\alpha, \gamma)$-contractive mapping, we get

$$\begin{aligned} d\left(x_{n(k)+1}, Tx^*\right) &= d\left(Tx_{n(k)}, Tx^*\right) \\ &\leq \alpha\left(x_{n(k)}, x^*\right) d\left(Tx_{n(k)}, Tx^*\right) \quad (20) \\ &\leq \gamma\left(d\left(x_{n(k)}, x^*\right)\right). \end{aligned}$$

Letting $k \rightarrow \infty$, we obtain $d(x_{n(k)+1}, Tx^*) \rightarrow 0$. Since $X$ is Hausdorff, we have that $Tx^* = x^*$.

This completes the proof. □

**Corollary 18.** *Let $(X, d)$ be a Hausdorff right $K$-complete $T_0$-quasi-pseudometric space. Suppose that $T : X \rightarrow X$ is an $(\alpha, \gamma)$-contractive mapping which satisfies the following conditions:*

(i) *$T$ is $\alpha$-admissible;*

(ii) *there exists $x_0 \in X$ such that $\alpha(Tx_0, x_0) \geq 1$;*

(iii) *if $(x_n)$ is a sequence in $X$ such that $\alpha(x_{n+1}, x_n) \geq 1$ for all $n = 1, 2, \ldots$ and $x_n \xrightarrow{d^{-1}} x$, then there exists a subsequence $(x_{n(k)})$ of $(x_n)$ such that $\alpha(x, x_{n(k)}) \geq 1$ for all $k$.*

*Then $T$ has a fixed point.*

**Corollary 19.** *Let $(X, d)$ be a bicomplete quasi-pseudometric space. Suppose that $T : X \rightarrow X$ is an $(\alpha, \gamma)$-contractive mapping which satisfies the following conditions:*

(i) *$T$ is $\alpha$-admissible and the function $\alpha$ is symmetric;*

(ii) *there exists $x_0 \in X$ such that $\alpha(Tx_0, x_0) \geq 1$;*

(iii) *if $(x_n)$ is a sequence in $X$ such that $\alpha(x_m, x_n) \geq 1$ for all $n, m \in \mathbb{N}$ and $x_n \xrightarrow{d^s} x$, then there exists a subsequence $(x_{n(k)})$ of $(x_n)$ such that $\alpha(x, x_{n(k)}) \geq 1$ for all $k$.*

*Then $T$ has a fixed point.*

## 4. Startpoint Theory

It is important to mention that there are a variety of endpoint concepts in the literature (see, e.g., [10]), each of them corresponding to a specified setting. Here we introduce a similar notion for set-valued maps defined on quasi-pseudometric spaces. Let $(X, d)$ be a $T_0$-quasi-pseudometric space.

*Definition 20.* Let $F : X \rightarrow 2^X$ be a set-valued map. An element $x \in X$ is said to be

(i) a fixed point of $F$ if $x \in Fx$,

(ii) a startpoint of $F$ if $H(\{x\}, Fx) = 0$,

(iii) an endpoint of $F$ if $H(Fx, \{x\}) = 0$,

(iv) an $\varepsilon$-startpoint of $F$ for some $\varepsilon \in (0, 1)$ if $H(\{x\}, Fx) < \varepsilon$,

(v) an $\varepsilon$-endpoint of $F$ for some $\varepsilon \in (0, 1)$ if $H(Fx, \{x\}) < \varepsilon$.

*Remark 21.* It is therefore obvious that if $x$ is both a startpoint of $F$ and an endpoint of $F$, then $x$ is a fixed point of $F$. In fact, $Fx$ is a singleton. But a fixed point need not be a startpoint nor an endpoint.

Indeed, consider the $T_0$-quasi-pseudometric space $(X, d)$, where $X = \{0, 1\}$ and $d$ is defined by $d(0, 1) = 0$, $d(1, 0) = 1$, and $d(x, x)$ for $x = 0, 1$. We define on $X$ the set-valued map $F : X \rightarrow 2^X$ by $Fx = X$. Obviously, 1 is a fixed point, but $H(\{1\}, F1) = H(\{1\}, X) = \max\{d(1, 1), d(1, 0)\} = 1 \neq 0$.

**Lemma 22.** *Let $(X, d)$ be a $T_0$-quasi-pseudometric space and let $F : X \to 2^X$ be a set-valued map. An element $x \in X$ is a startpoint of $F$ if and only if it is an $\varepsilon$-startpoint of $F$ for every $\varepsilon \in (0, 1)$.*

**Lemma 23.** *Let $(X, d)$ be a $T_0$-quasi-pseudometric space and let $F : X \to 2^X$ be a set-valued map. An element $x \in X$ is an endpoint of $F$ if and only if it is an $\varepsilon$-endpoint of $F$ for every $\varepsilon \in (0, 1)$.*

*Definition 24.* Let $(X, d)$ be a $T_0$-quasi-pseudometric space. We say that a set-valued map $F : X \to 2^X$ has the *approximate startpoint property* (resp., *approximate endpoint property*) if

$$\inf_{x \in X} \sup_{y \in Fx} d(x, y) = 0 \quad \left(\text{resp.}, \inf_{x \in X} \sup_{y \in Fx} d(y, x) = 0\right). \quad (21)$$

*Definition 25.* Let $(X, d)$ be a $T_0$-quasi-pseudometric space. We say that a set-valued map $F : X \to 2^X$ has the *approximate mix-point property* if

$$\inf_{x \in X} \sup_{y \in Fx} d^s(x, y) = 0. \quad (22)$$

Here, it is also very clear that $F$ has approximate mix-point property if and only if $F$ has both the approximate startpoint and the approximate endpoint properties.

We are therefore naturally led to this definition.

*Definition 26.* Let $f : X \to X$ be a single-valued map on a $T_0$-quasi-pseudometric space $(X, d)$. Then $f$ has the *approximate startpoint property* (resp., *approximate endpoint property*) if and only if

$$\inf_{x \in X} d(x, fx) = 0 \quad \left(\text{resp.}, \inf_{x \in X} d(fx, x) = 0\right). \quad (23)$$

We motivate our coming results by the following examples. We basically show that the concepts defined above are independent and do not necessarily coincide. The list of examples presented is not exhaustive and more can be constructed, showing the connection between the notions defined above.

*Example 27.* Let $X = \{0, 1, 2\}$. The map $d : X \times X \to [0, \infty)$ defined by $d(0, 1) = d(0, 2) = 0$, $d(1, 0) = d(1, 2) = 1$, $d(2, 0) = d(2, 1) = 2$, and $d(x, x) = 0$ for all $x \in X$ is a $T_0$-quasi-pseudometric on $X$. Let $F : X \to 2^X$ be the set mapping defined by $Fa = X \setminus \{a\}$ for any $a \in X$. By definition, $F$ does not have any fixed point. Nevertheless, a simple computation shows that $H(\{0\}, F0) = 0$, and hence $0$ is a startpoint and it is the only one. Also there is no endpoint. Again, with a direct computation, we have $\inf_{x \in X} \sup_{y \in Fx} d(x, y) = 0$, showing that $F$ has the approximate startpoint property, but $\inf_{x \in X} \sup_{y \in Fx} d(y, x) = 1$, showing that $F$ does not have the approximate endpoint property.

*Example 28.* Let $X = \{1/n, n = 1, 2, \ldots\}$. The map $d : X \times X \to [0, \infty)$ defined by $d(1/n, 1/m) = \max\{1/n - 1/m, 0\}$ is a $T_0$-quasi-pseudometric on $X$. Let $F : X \to 2^X$ be the set-valued mapping defined by $Fa = X \setminus \{a\}$ for any $a \in X$. By definition, $F$ does not have any fixed point.

For a fixed $n_0 \in \mathbb{N}$,

$$d\left(\frac{1}{n_0}, \frac{1}{m}\right) = \begin{cases} 0 & \text{if } m \leq n_0, \\ \dfrac{1}{n_0} - \dfrac{1}{m} & \text{if } m > n_0. \end{cases} \quad (24)$$

Similarly for a fixed $n_0 \in \mathbb{N}$,

$$d\left(\frac{1}{m}, \frac{1}{n_0}\right) = \begin{cases} 0 & \text{if } m \geq n_0, \\ \dfrac{1}{m} - \dfrac{1}{n_0} & \text{if } m < n_0. \end{cases} \quad (25)$$

Hence, $F$ does not have any startpoint nor endpoint (which also implies that $F$ does not have any fixed point).

But for a given $\varepsilon \in (0, 1)$, there exists $k_0 \in \mathbb{N}$ such that $1/k_0 < \varepsilon$. We also know from definition that $H(\{1/k_0\}, F(1/k_0)) = 1/k_0$, so, $1/k_0$ is an $\varepsilon$-startpoint of $F$. More generally, for a given $\varepsilon \in (0, 1)$, there exists $k_\varepsilon \in \mathbb{N}$ such that $1/k_\varepsilon$ is an $\varepsilon$-startpoint of $F$. Moreover, for any $k \geq k_\varepsilon$, $1/k$ is an $\varepsilon$-startpoint of $F$.

Similarly, we can show that $F$ admits an $\varepsilon$-endpoint.

We can now state our first result.

**Theorem 29.** *Let $(X, d)$ be a bicomplete quasi-pseudometric space. Let $F : X \to CB(X)$ be a set-valued map that satisfies*

$$H(Fx, Fy) \leq \psi(d(x, y)), \quad \text{for each } x, y \in X, \quad (26)$$

*where $\psi : [0, \infty) \to [0, \infty)$ is upper semicontinuous, $\psi(t) < t$ for each $t > 0$, and $\lim \inf_{t \to \infty} (t - \psi(t)) > 0$. Then there exists a unique $x_0 \in X$ which is both a startpoint and an endpoint of $F$ if and only if $F$ has the approximate mix-point property.*

*Proof.* It is clear that if $F$ admits a point which is both a startpoint and an endpoint, then $F$ has the approximate startpoint property and the approximate endpoint property. Just observe that $H(\{x\}, Fx) = \sup_{y \in Fx} d(x, y)$ and $H(Fx, \{x\}) = \sup_{y \in Fx} d(y, x)$. Conversely, suppose $F$ has the approximate mix-point property. Then

$$C_n = \left\{ x \in X : \sup_{y \in Fx} d^s(x, y) \leq \frac{1}{n} \right\} \neq \emptyset, \quad (27)$$

for each $n \in \mathbb{N}$. Also it is clear that for each $n \in \mathbb{N}$, $C_{n+1} \subseteq C_n$. The map $x \mapsto \sup_{y \in Fx} d^s(x, y)$ is $\tau(d^s)$-lower semicontinuous (as supremum of $\tau(d^s)$-continuous mappings); we have that $C_n$ is $\tau(d^s)$-closed.

Next we prove that, for each $n \in \mathbb{N}$, $C_n$ is bounded.

Assume by the way of contradiction that $\delta(C_n) = \infty$ for each $n \in \mathbb{N}$. Then there exist $x_n, y_n \in C_n$ such that $d(x_n, y_n) \geq n$. From (26), we obtain that

$$
\begin{aligned}
d(x, y) &= H(\{x\}, \{y\}) \\
&\leq H(\{x\}, Fx) + H(Fx, Fy) + H(Fy, \{y\}) \\
&= \frac{2}{n} + \psi(d(x, y)),
\end{aligned}
\tag{28}
$$

whenever $x, y \in C_n$.

Therefore,

$$
d(x, y) - \psi(d(x, y)) \leq \frac{2}{n},
\tag{29}
$$

whenever $x, y \in C_n$. Hence,

$$
\begin{aligned}
\lim_{n \to \infty} (d(x_n, y_n) - \psi(d(x_n, y_n))) &= 0, \\
\lim_{n \to \infty} d(x_n, y_n) &= \infty.
\end{aligned}
\tag{30}
$$

This contradicts our assumption. Now we show that $\lim_{n \to \infty} \delta(C_n) = 0$. On the contrary, assume $\lim_{n \to \infty} \delta(C_n) = r_0 > 0$ (note that the sequence $(\delta(C_n))$ is nonincreasing and bounded below and then has a limit). Let

$$
\begin{aligned}
\beta = \inf_{n \in \mathbb{N}} \inf \Big\{ &\liminf_{k \to \infty} (r_{n,k} - \psi(r_{n,k})) : (x_{n,k}, y_{n,k}) \in C_n, \\
&r_{n,k} = d(x_{n,k}, y_{n,k}) \longrightarrow \delta(C_n) \text{ as } k \longrightarrow \infty \Big\}.
\end{aligned}
\tag{31}
$$

Now we show that $\beta > 0$ (notice $\beta \geq 0$). Arguing by contradiction, we assume $\beta = 0$; then by the definition of $\beta$, there exists a sequence $r_n$ such that $r_n \to r_0$ and $\lim_{n \to \infty} (r_n - \psi(r_n)) = 0$. Then $\lim_{n \to \infty} \psi(r_n) = r_0$. But since $\psi$ is upper semicontinuous and $r_0 > 0$, then

$$
r_0 = \lim_{n \to \infty} \psi(r_n) \leq \psi(r_0) < r_0.
\tag{32}
$$

This contradiction shows that $\beta > 0$. For each $n \in N$, let $(x_k, y_k) \in C_n$ be a sequence such that $d(x_k, y_k) \to \delta(C_n)$, as $k \to \infty$. Then from (29) we get

$$
\beta \leq \liminf_{k \to \infty} (d(x_k, y_k) - \psi(d(x_k, y_k))) \leq \frac{2}{n},
\tag{33}
$$
$$
\text{for each } n \in N.
$$

Hence, $\beta = 0$. This contradiction shows that $\lim_{n \to \infty} \delta(C_n) = 0$. It follows from the Cantor intersection theorem that $\bigcap_{n \in N} C_n = \{x_0\}$.

Thus, $H(\{x_0\}, Fx_0) = \sup_{y \in Fx_0} d(x, y) = 0 = \sup_{y \in Fx_0} d(y, x) = H(Fx_0, \{x_0\})$. For uniqueness, if $z_0$ is an arbitrary startpoint of $F$, then $H(\{z_0\}, Fz_0) = 0 = H(Fz_0, \{z_0\})$, and so $z_0 \in \bigcap_{n \in N} C_n = \{x_0\}$.

This completes the proof. □

**Corollary 30.** *Let $(X, d)$ be a bicomplete quasi-pseudometric space. Let $F : X \to CB(X)$ be a set-valued map that satisfies*

$$
H(Fx, Fy) \leq \psi(d(x, y)), \quad \text{for each } x, y \in X,
\tag{34}
$$

*where $\psi : [0, \infty) \to [0, \infty)$ is upper semicontinuous, $\psi(t) < t$ for each $t > 0$, and $\liminf_{t \to \infty} (t - \psi(t)) > 0$. If $F$ has the approximate mix-point property then $F$ has a fixed point.*

*Proof.* From Theorem 29, we conclude that there exists $x_0$ which is both a startpoint and an endpoint; that is, $H(\{x_0\}, Fx_0) = 0 = H(Fx_0, \{x_0\})$. The $T_0$-condition therefore guarantees the desired result. □

**Theorem 31.** *Let $(X, d)$ be a bicomplete quasi-pseudometric space. Let $F : X \to CB(X)$ be a set-valued map that satisfies*

$$
H(Fx, Fy) \leq kd(x, y), \quad \text{for each } x, y \in X,
\tag{35}
$$

*where $0 \leq k < 1$. Then there exists a unique $x_0 \in X$ which is both a startpoint and an endpoint of $F$ if and only if $F$ has the approximate mix-point property.*

*Proof.* Take $\psi(t) = kt$ in Theorem 29. □

We then deduce the following result for single-valued maps.

**Theorem 32.** *Let $(X, d)$ be a bicomplete quasi-pseudometric space. Let $f : X \to X$ be a map that satisfies*

$$
d(fx, fy) \leq \psi(d(x, y)), \quad \text{for each } x, y \in X,
\tag{36}
$$

*where $\psi : [0, \infty) \to [0, \infty)$ is upper semicontinuous, $\psi(t) < t$ for each $t > 0$, and $\liminf_{t \to \infty} (t - \psi(t)) > 0$. Then $f$ has the approximate startpoint property.*

*Proof.* By the way of contradiction, suppose that $\inf_{x \in X} d(x, fx) > 0$. Then

$$
\begin{aligned}
\inf_{x \in X} d(x, fx) &\leq \inf_{y \in f(X)} d(y, fy) \\
&= \inf_{x \in X} d(fx, f^2 x) \\
&\leq \inf_{x \in X} \psi(d(x, fx)).
\end{aligned}
\tag{37}
$$

Since $\psi(d(x, y)) \leq d(x, y)$, then $\inf_{x \in X} \psi(d(x, fx)) \leq \inf_{x \in X} d(x, fx)$.

Now, on the contrary, suppose again that

$$
\inf_{x \in X} \psi(d(x, fx)) = \inf_{x \in X} d(x, fx).
\tag{38}
$$

Let $(x_n) \subset X$ be a sequence such that $\lim_{n \to \infty} d(x_n, f(x_n)) = \inf_{x \in X} d(x, fx)$. By passing to subsequences if necessary, we may assume that $\lim_{n\to\infty}\psi(d(x_n, f(x_n)))$ exists. Then from (37) we have

$$
\begin{aligned}
\inf_{x\in X} d(x, fx) &\leq \inf_{x\in X}\psi(d(x, fx)) \\
&\leq \lim_{n\to\infty}\psi(d(x_n, f(x_n))) \\
&\leq \psi\left(\inf_{x\in X} d(x, fx)\right) \\
&< \inf_{x\in X} d(x, fx).
\end{aligned}
\tag{39}
$$

We get a contradiction, so $\inf_{x\in X}\psi(d(x, fx)) < \inf_{x\in X} d(x, fx)$ which again contradicts (37).

This completes the proof. □

**Corollary 33.** *Let $(X, d)$ be a bicomplete quasi-pseudometric space. Let $f : X \to X$ be a map that satisfies*

$$
d(fx, fy) \leq \psi(d(x, y)), \quad \text{for each } x, y \in X,
\tag{40}
$$

*where $\psi : [0, \infty) \to [0, \infty)$ is upper semicontinuous, $\psi(t) < t$ for each $t > 0$, and $\liminf_{t\to\infty}(t - \psi(t)) > 0$. Then $f$ has the approximate endpoint property.*

We finish this section by the following fixed point result.

**Corollary 34.** *Let $(X, d)$ be a bicomplete quasi-pseudometric space. Let $f : X \to X$ be a map that satisfies*

$$
d(fx, fy) \leq \psi(d(x, y)), \quad \text{for each } x, y \in X,
\tag{41}
$$

*where $\psi : [0, \infty) \to [0, \infty)$ is upper semicontinuous, $\psi(t) < t$ for each $t > 0$, and $\liminf_{t\to\infty}(t - \psi(t)) > 0$. Then $f$ has a fixed point.*

*Proof.* From Theorem 32 and Corollary 33, we conclude that $f$ has the approximate mix-point property. Hence, by Corollary 30, we have the desired result. □

## 5. More Results

The following theorem is the main result of this section.

**Theorem 35.** *Let $(X, d)$ be a left $K$-complete quasi-pseudometric space. Let $F : X \to CB(X)$ be a set-valued map and $f : X \to \mathbb{R}$ as $f(x) = H(\{x\}, Fx)$. If there exists $c \in (0, 1)$ such that for all $x \in X$ there exists $y \in Fx$ satisfying*

$$
H(\{y\}, Fy) \leq c(d(x, y)),
\tag{42}
$$

*then $T$ has a startpoint.*

*Proof.* For any initial $x_0 \in X$, there exists $x_1 \in Fx_0 \subseteq X$ such that

$$
H(\{x_1\}, Fx_1) \leq c(d(x_0, x_1)),
\tag{43}
$$

and for $x_1 \in X$, there is $x_2 \in Fx_1 \subseteq X$ such that

$$
H(\{x_2\}, Fx_2) \leq c(d(x_1, x_2)).
\tag{44}
$$

Continuing this process, we can get an iterative sequence $(x_n)$ where $x_{n+1} \in Fx_n \subseteq X$ and

$$
H(\{x_{n+1}\}, Fx_{n+1}) \leq c(d(x_n, x_{n+1})), \quad n = 0, 1, 2, \ldots.
\tag{45}
$$

*Claim 1.* $(x_n)$ is a left $K$-Cauchy sequence.

On one hand,

$$
H(\{x_{n+1}\}, Fx_{n+1}) \leq c(d(x_n, x_{n+1})), \quad n = 0, 1, 2, \ldots.
\tag{46}
$$

On the other hand, $x_{n+1} \in Fx_n$ implies

$$
d(x_n, x_{n+1}) \leq H(\{x_n\}, Fx_n), \quad n = 0, 1, 2, \ldots.
\tag{47}
$$

By the two above inequalities, we have

$$
d(x_{n+1}, x_{n+2}) \leq cd(x_n, x_{n+1}) \quad n = 0, 1, 2, \ldots,
$$
$$
H(\{x_{n+1}\}, Fx_{n+1}) \leq cH(\{x_n\}, Fx_n) \quad n = 0, 1, 2, \ldots.
\tag{48}
$$

We then deduce by iteration that

$$
d(x_n, x_{n+1}) \leq c^n d(x_0, x_1) \quad n = 0, 1, 2, \ldots,
$$
$$
H(\{x_{n+1}\}, Fx_{n+1}) \leq c^n H(\{x_0\}, Fx_0) \quad n = 0, 1, 2, \ldots.
\tag{49}
$$

Then for $n, m \in \mathbb{N}, n < m$,

$$
\begin{aligned}
d(x_n, x_m) &\leq d(x_n, x_{n+1}) + d(x_{n+1}, x_{n+2}) \\
&\quad + \cdots + d(x_{m-1}, x_m) \\
&\leq \left[c^n + c^{n+1} + \cdots + c^{m-1}\right] d(x_0, x_1) \\
&\leq \frac{c^n}{1-c} d(x_0, x_1),
\end{aligned}
\tag{50}
$$

and since $c^n \to 0$ as $n \to \infty$ we conclude that $(x_n)$ is a left $K$-Cauchy sequence. According to the left $K$-completeness of $(X, d)$, there exists $x^* \in X$ such that $x_n \xrightarrow{d} x^*$.

*Claim 2.* $x^*$ is a startpoint of $F$.

Observe that the sequence $(fx_n) = (H(\{x_n\}, Fx_n))$ is decreasing and hence converges to 0. Since $f$ is $\tau(d)$-lower semicontinuous (as supremum of $\tau(d)$-lower semicontinuous functions), we have

$$
0 \leq f(x^*) \leq \liminf_{n\to\infty} f(x_n) = 0.
\tag{51}
$$

Hence, $f(x^*) = 0$; that is, $H(\{x^*\}, Fx^*) = 0$.

This completes the proof. □

*Example 36.* Let $X = \{1, 1/2, 1/3\}$. The map $d : X \times X \to [0, \infty)$ defined by $d(1/n, 1/m) = \max\{1/n - 1/m, 0\}$ is a $T_0$-quasi-pseudometric on $X$. Let $F : X \to 2^X$ be the set mapping defined by $Fa = X \setminus \{a\}$ for any $a \in X$. With $c = 1/2$, the map $F$ satisfies the assumptions of our theorem, so it has a startpoint, which in the present case is $1/3$.

More generally, if we set $X_n = \{1/i, \ i = 1, 2, \ldots, n\}$ and $d$ as defined above, with $c = 1/2$, the map $F$ defined by $Fa = X \setminus \{a\}$ for any $a \in X$ satisfies the assumptions of our theorem, so it has a startpoint, which in this case is $1/n$.

**Corollary 37.** *Let $(X, d)$ be a right $K$-complete quasi-pseudometric space. Let $F : X \rightarrow CB(X)$ be a set-valued map and $f : X \rightarrow \mathbb{R}$ defined by $f(x) = H(Fx, \{x\})$. If there exists $c \in (0, 1)$ such that for all $x \in X$ there exists $y \in Fx$ satisfying*

$$H\left(Fy, \{y\}\right) \le c\left(d\left(y, x\right)\right), \tag{52}$$

*then $T$ has an endpoint.*

**Corollary 38.** *Let $(X, d)$ be a bicomplete quasi-pseudometric space. Let $F : X \rightarrow CB(X)$ be a set-valued map and $f : X \rightarrow \mathbb{R}$ defined by $f(x) = H^s(Fx, \{x\}) = \max\{H(Fx, \{x\}), H(\{x\}, Fx)\}$. If there exists $c \in (0, 1)$ such that for all $x \in X$ there exists $y \in Fx$ satisfying*

$$H^s\left(\{y\}, Fy\right) \le c\left(\min\left\{d\left(y, x\right), d\left(y, x\right)\right\}\right), \tag{53}$$

*then $T$ has a fixed point.*

*Proof.* We give here the main idea of the proof.

Observe that inequality (53) guarantees that the sequence $(x_n)$ constructed in the proof of Theorem 35 is a $d^s$-Cauchy sequence and hence $d^s$-converges to some $x^*$. Using the fact that $f$ is $\tau(d^s)$-lower semicontinuous (as supremum of $\tau(d^s)$-continuous functions), we have

$$0 \le f\left(x^*\right) \le \liminf_{n \to \infty} f\left(x_n\right) = 0. \tag{54}$$

Hence, $f(x^*) = 0$; that is, $H(\{x^*\}, Fx^*) = 0 = H(Fx^*, \{x^*\})$, and we are done. $\qquad\square$

*Remark 39.* All the results given remain true when we replace accordingly the bicomplete quasi-pseudometric space $(X, d)$ with a left Smyth sequentially complete/left $K$-complete or a right Smyth sequentially complete/right $K$-complete space.

## Conflict of Interests

The author declares that there is no conflict of interests regarding the publication of this paper.

## References

[1] K. Włodarczyk and R. Plebaniak, "Asymmetric structures, discontinuous contractions and iterative approximation of fixed and periodic points," *Fixed Point Theory and Applications*, vol. 2013, article 128, 2013.

[2] H.-P. A. Künzi, "An introduction to quasi-uniform spaces," *Contemporary Mathematics*, vol. 486, pp. 239–304, 2009.

[3] K. Wlodarczyk and R. Plebaniak, "Generalized uniform spaces, uniformly locally contractive set-valued dynamic systems and fixed points," *Fixed Point Theory and Applications*, vol. 2012, article 104, 2012.

[4] K. W. lodarczyk and R. Plebaniak, "Fixed points and endpoints of contractive set-valued maps in cone uniform spaces with generalized pseudo distances," *Fixed Point Theory and Applications*, vol. 2012, article 176, 15 pages, 2012.

[5] K. W. Włodarczyk and R. Plebaniak, "Leader type contractions, periodic and fixed points and new completivity in quasi-gauge spaces with generalized quasi-pseudodistances," *Topology and its Applications*, vol. 159, no. 16, pp. 3504–3512, 2012.

[6] K. Włodarczyk and R. Plebaniak, "Contractivity of Leader type and fixed points in uniform spaces with generalized pseudodistances," *Journal of Mathematical Analysis and Applications*, vol. 387, no. 2, pp. 533–541, 2012.

[7] K. Włodarczyk and R. Plebaniak, "Contractions of Banach, Tarafdar, Meir-Keeler, Ćirić-Jachymski-Matkowski and Suzuki types and fixed points in uniform spaces with generalized pseudodistances," *Journal of Mathematical Analysis and Applications*, vol. 404, no. 2, pp. 338–350, 2013.

[8] K. W lodarczyk and R. Plebaniak, "New completeness and periodic points of discontinuous contractions of Banach type in quasi-gauge spaces without Hausdorff property," *Fixed Point Theory and Applications*, vol. 2013, article 289, 27 pages, 2013.

[9] H.-P. Künzi and C. Ryser, "The Bourbaki quasi-uniformity," *Topology Proceedings*, vol. 20, pp. 161–183, 1995.

[10] C. A. Agyingi, P. Haihambo, and H.-P. A. Künzi, "Endpoints in $T_0$-quasimetric spaces: part II," *Abstract and Applied Analysis*, vol. 2013, Article ID 539573, 10 pages, 2013.

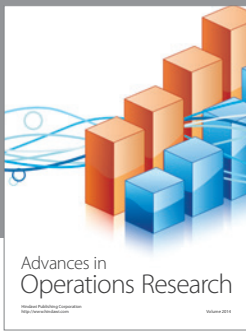
Advances in
Operations Research

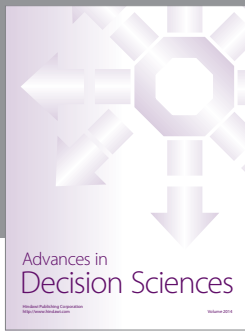
Advances in
Decision Sciences

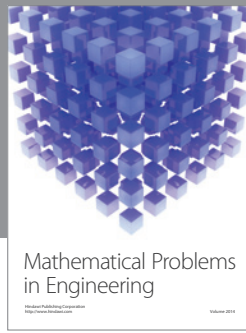
Mathematical Problems
in Engineering

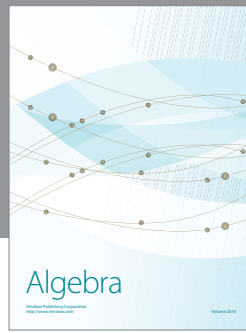
Algebra

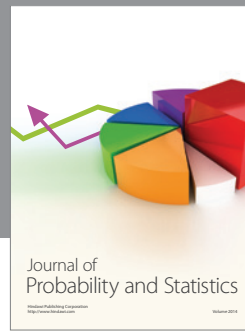
Journal of
Probability and Statistics

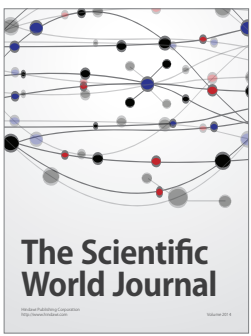
The Scientific
World Journal

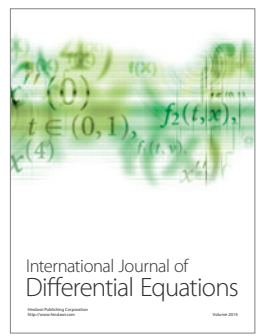
International Journal of
Differential Equations

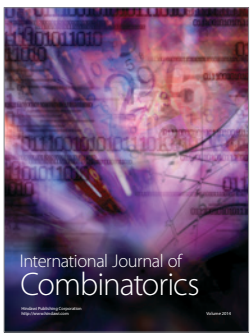
International Journal of
Combinatorics

Submit your manuscripts at
http://www.hindawi.com

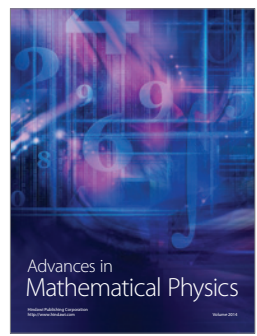
Advances in
Mathematical Physics

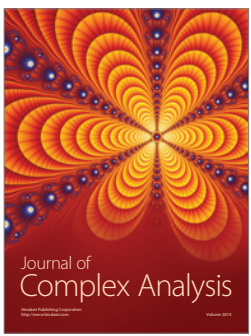
Journal of
Complex Analysis

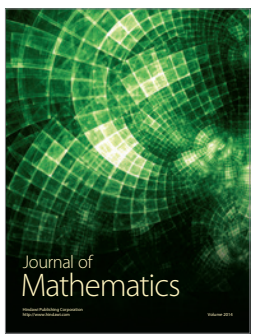
Journal of
Mathematics

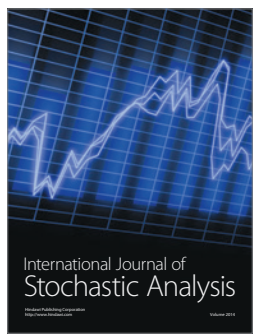
International Journal of
Stochastic Analysis

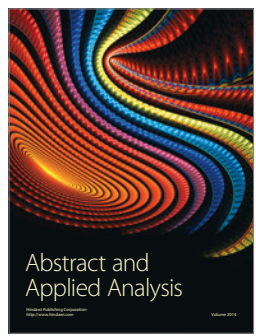
Abstract and
Applied Analysis

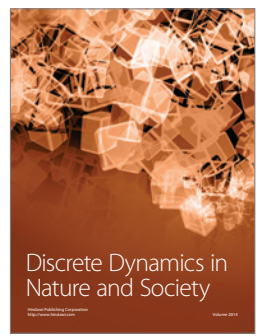
Discrete Dynamics in
Nature and Society

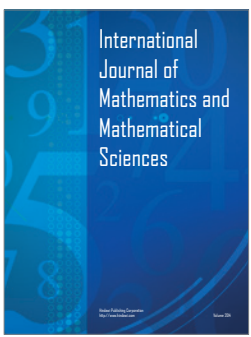
International
Journal of
Mathematics and
Mathematical
Sciences

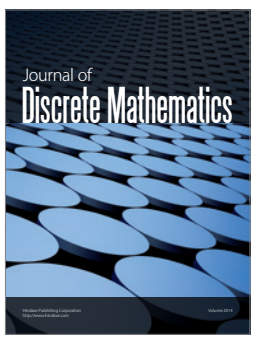
Journal of
Discrete Mathematics

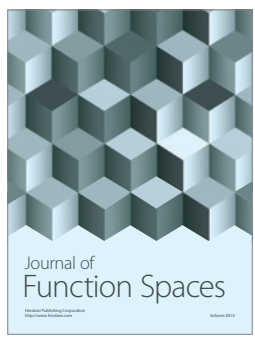
Journal of
Function Spaces

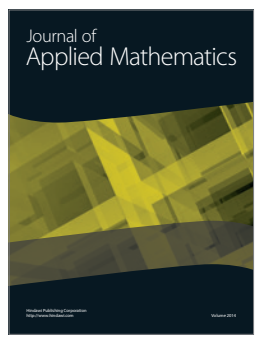
Journal of
Applied Mathematics

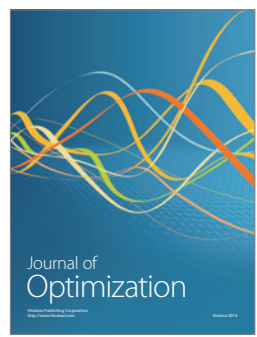
Journal of
Optimization